\documentclass[11pt,reqno,oneside]{amsart}
\usepackage{amscd}
\usepackage{amsmath}
\usepackage{latexsym}
\usepackage{amsfonts}
\usepackage{amssymb}
\usepackage{amsthm}
\usepackage{graphicx}
\usepackage{verbatim}
\usepackage{mathrsfs}
\usepackage{enumerate}
\usepackage{hyperref}
\usepackage{fullpage}
\usepackage{xcolor}
\usepackage[latin1]{inputenc}
 
\theoremstyle{plain}
\theoremstyle{plain}\newtheorem{theorem}{Theorem}[section]
\theoremstyle{plain}\newtheorem{lemma}[theorem]{Lemma}
\theoremstyle{plain}
\theoremstyle{plain}
\theoremstyle{plain}\newtheorem{remark}{Remark}[section]

\numberwithin{equation}{section}

\newcommand{\R}{\mathbb{R}}
\newcommand{\be}{\begin{equation}}
\newcommand{\ee}{\end{equation}}
 \newcommand{\ba}{\begin{aligned}}
 \newcommand{\ea}{\end{aligned}}

  \newcommand{\ben}{\begin{enumerate}}
   \newcommand{\een}{\end{enumerate}}

\newcommand{\Rmnum}[1]{\expandafter\@slowromancap\romannumeral #1@}
\allowdisplaybreaks

\begin{document}

\title{Formation of singularities for  multi-dimensional transport equations with nonlocal velocity}
\author[]{Quansen Jiu$^{1}$,\, Wanwan Zhang$^{2}$}

\address{$^1$ School of Mathematical Sciences, Capital Normal University, Beijing, 100048, P.R.China}

\email{jiuqs@cnu.edu.cn}

\address{$^{2}$School of Mathematical Sciences, Capital Normal University, Beijing, 100048, P.R.China}
\email{zhangwanwan153@163.com}

\subjclass[2000]{35Q35; 35B35; 76D05}
\keywords{Transport equations; Nonlocal velocity; Formation of singularities; Riesz transform}

\begin{abstract}
This paper is concerned with a class of multi-dimensional transport equations with nonlocal velocity. It is shown that the local smooth solution cannot exist globally in time via the De Giorgi iteration technique.
\end{abstract}
\smallskip
\maketitle


\section{Introduction}
In this paper, we consider the Cauchy problem of the following multi-dimensional transport equation with nonlocal velocity
\begin{equation}\label{1.1}
\left\{\ba
&\theta_{t}+u\cdot\nabla \theta = 0, ~(x,t)\in \R^{n}\times\R^{+},\\&u=\nabla\Lambda^{-2+2\alpha}\theta,\\
&\theta(x,0)=\theta_{0}(x),~x\in\R^n, \ea\ \right.
\end{equation}
with $0<\alpha<1$ and $n\geq2$. Here the unknown function is $\theta=\theta(x,t)$ with $x=(x_1,...,x_n)\in\R^n$, $t>0$, $\nabla=(\partial_{x_1},...,\partial_{x_n})$ and the fractional Laplacian $\Lambda^s=(-\Delta)^{\frac{s}{2}}$ with $s\in \R$ is defined through the Fourier transform
\begin{eqnarray} \label{1.2}
\widehat{\Lambda^sf}(\xi)=|\xi|^s\widehat{f}(\xi),
\end{eqnarray}
where
\begin{eqnarray*}
\widehat{f}(\xi)=\int_{\R^n}e^{-2\pi ix\cdot\xi}f(x)\,dx.
\end{eqnarray*}
The corresponding dissipative equation reads as
\begin{equation}\label{1.3}
\left\{\ba
&\theta_{t}+u\cdot\nabla \theta+\kappa\Lambda^\gamma\theta = 0, ~(x,t)\in \R^{n}\times\R^{+}\\&u=\nabla\Lambda^{-2+2\alpha}\theta\\
&\theta(x,0)=\theta_{0}(x),~x\in\R^{n} \ea\ \right.
\end{equation}
with $0<\alpha<1$, $0<\gamma<2$ and $\kappa>0$.

It is known that the velocity field $u$ in the second equation of \eqref{1.1} can  be expressed as  (see \cite{[S]}):
\begin{eqnarray}\label{1.4}
u(x,t)=C_{n,\alpha}P.V.\int_{\R^{n}}\frac{x-y}{|x-y|^{n+2\alpha}}\theta(y,t)dy,
\end{eqnarray}
where
$C_{n,\alpha}=\frac{(2-2\alpha-n)\Gamma(\frac{n}{2}-1+\alpha)}{\pi^{\frac{n}{2}}2^{2-2\alpha}\Gamma(1-\alpha)}$
with  $\Gamma$  the  Gamma function, and P.V. denotes the principle value integration.

The transport equation  \eqref{1.1} is a class of nonlocal active scalar equation. The question of global regularity or finite time singularity for active scalar equation with nonlocal velocity field has attracted much attention in recent years. We refer the reader to a remarkable review paper \cite{[K]} and the references therein for  recent progress in this area.
One motivation of \eqref{1.1} is that it can be regarded as a model equation for understanding the generalized surface quasi-geostrophic (GSQG) equation, in which the velocity field $u$ is given by
\begin{eqnarray}\label{1.5}
u(x,t)=\nabla^\bot\Lambda^{-2+2\alpha}\theta(x,t),
\end{eqnarray}
where $\nabla^\bot=(\partial_{x_2},-\partial_{x_1})$. Clearly, the velocity \eqref{1.5} in the GSQG equation is divergence free, while the velocity of \eqref{1.1} is not in general.
In the case $\alpha=\frac{1}{2}$ in \eqref{1.5}, GSQG equation is reduced to the classical surface quasi-geostrophic equation (SQG). There are a number of  mathematical studies on GSQG and SQG equations and we refer the readers  to \cite{[LCV],[CMT],[CTV],[CV],[CW],[CC],[C],[CotVv],[KNV],[KRYZ],[KYZ],[MB],[P]} and \cite{[R]} for more details. Here we briefly recall some recent results for \eqref{1.1} and \eqref{1.3}.
Dong and Li \cite{[DL]} proved that certain radial solutions must develop gradient blow-up in finite time for  \eqref{1.1} with $\alpha=\frac{1}{2}$.
 In \cite{[D]}, Dong was able to obtain the blow-up of smooth radial solutions to \eqref{1.1} with more general $\alpha\in (0,1)$ for any smooth, radially symmetric and nonnegative initial data with compact support and its positive maximum at the origin.  Moreover, Li and Rodrigo \cite{[LR09]} proved formation of singularities of solutions to \eqref{1.3} with $\frac{1}{4}\leq\alpha\leq\frac{1}{2}$ and $0\leq\gamma<\alpha$ for a generic family of initial data. It is  conjectured in \cite{[LR09]} that singularities exist for \eqref{1.3} with $\alpha=\frac{1}{2}$ and $\frac{1}{2}\leq\gamma<1$, which remains  open.

Another example of active scalar equation with nonlocal velocity field is the following CCF model
\begin{eqnarray}\label{1.6}
\partial_t\theta+H(\theta)\theta_x+\kappa\Lambda^\gamma\theta=0,
\end{eqnarray}
where $H$ is the Hilbert transform defined by
\begin{eqnarray*}
H(\theta)(x)=\frac{1}{\pi}P.V.\int_{\R}\frac{\theta(y)}{x-y}dy.
\end{eqnarray*}
This model was introduced by C\'{o}rdoba, C\'{o}rdoba and Fontelos in the pioneering work \cite{[CCF05]}.
By virtue of the Meillin transform and the complex analysis, they obtained some new bilinear
estimates for the Hilbert transform and as a result proved the breakdown of local smooth solutions to \eqref{1.6}
with $\kappa=0$ for a generic class of initial data.
In \cite{[CCF06]}, more and thorough integral inequalities involving the Hilbert transform are studied, from which they obtained the finite time singularity of solutions to the CCF model without viscosity.
Later, Silvestre and Vicol \cite{[SV]}  provided four  different and interesting proofs of the finite time singularity formulation of solution to \eqref{1.6} with $\kappa=0$, which was first proved in \cite{[CCF05]}.
In \cite{[LR08]}, Li and Rodrigo were able to extend the inviscid blow up result to include the dissipative term
for $0\leq\gamma<\frac{1}{2}$. Another proof of this dissipative singularity result via the telescopic sum argument is presented in \cite{[SV]}. Recently, in \cite{[Lucas]}, Ferreira and Moitinho obtained the existence of global classical solutions to \eqref{1.6} with $\gamma\in(\gamma_1,1)$, where $\gamma_1$ depends on some norm of the initial data.
The question of global regularity vs finite time singularity for the dissipative CCF model \eqref{1.6} with $\frac{1}{2}\leq\gamma<1$ is  open as well.
For more progress in this direction, we refer the reader to \cite{[LR20]}. The multi-dimensional extension on the CCF model  \eqref{1.6}
with $\kappa=0$  was done by Balodis and C\'{o}rdoba in \cite{[BC]}, in which the authors investigated the equation
\begin{eqnarray}\label{1.7}
\partial_t\theta+\mathcal{R}(\theta)\cdot\nabla\theta=0.
\end{eqnarray}
  Here the Riesz transform $\mathcal{R}$ is defined by
\begin{eqnarray*}
\mathcal{R}=(\mathcal{R}_1,\mathcal{R}_2,...,\mathcal{R}_n)=(\frac{\partial_{x_1}}{\sqrt{-\Delta}},\frac{\partial_{x_2}}{\sqrt{-\Delta}},...,\frac{\partial_{x_n}}{\sqrt{-\Delta}}),
\end{eqnarray*}
which can  be written with the aid of the Fourier transform as $\widehat{\mathcal{R}f}(\xi)=i\frac{\xi}{|\xi|}\widehat{f}(\xi)$.
In \cite{[BC]}, the authors proved the local well-posedness to \eqref{1.7} in Sobolev space $H^s(\R^n)$ with $s>1+\frac{n}{2}$ and  obtained the blow-up of the local smooth solution for any nonnegative, not-identically zero initial data. The model \eqref{1.7} can be also seen as  one of \eqref{1.1} with $\alpha=\frac{1}{2}$.

In this paper, motivated by \cite{[SV]}, we will prove the finite time singularity of the local smooth solution to \eqref{1.1} with full range $0<\alpha<1$ for more general initial data via the De Giorgi iteration technique. Our main result can be stated as
\begin{theorem}\label{the-1}
For any initial data $\theta_0\in \mathcal{S}(\R^n)$, the Schwartz function class, satisfying $\displaystyle\sup_{x\in\R^n}\theta_0(x)>0$, there is no smooth global solution $\theta$ to \eqref{1.1}.
\end{theorem}
Let us describe the technique in some detail. It is clear that the integral of solutions to GSQG equations is a conservative quantity. However, \eqref{1.1} is not an incompressible model in general.
One of main idea is to make use of the  structure of nonlinearity term, which   makes the $L^1(\R^n)$ norm  of the nonnegative solution dissipative (see Lemma \ref{identity}).
More precisely, one can obtain the following inequality concerning the level of truncations $\theta_k$ (see \eqref{4.1} for the definition of $\theta_k$)
$$\frac{d}{dt}\|\theta_k(\cdot,t)\|_{L^1(\R^n)}\leq-\|\theta_k(\cdot,t)\|^2_{\dot{H}^\alpha(\R^n)},$$
which  plays a similar role as an energy dissipation inequality in the De Giorgi iteration scheme. Furthermore, by virtue of the De Giorgi iteration technique, we are able to obtain a decay for
$\displaystyle\sup_{x\in\mathbb{R}^{n}}\theta(x,t)$, which implies that smooth solutions must blow up in finite time.
In comparison with  \cite{[D]}, the initial data  is not necessarily radial symmetric, nonnegative and with its positive maximum at the origin in Theorem \ref{the-1} and  a  different approach is given in this paper.

The remaining part of this paper is organized as follows.  In Section 2 we recall the local well-posedness result of \eqref{1.1}, a recurrence lemma and some interpolation inequality. Section 3 is devoted to some basic properties of the solutions to \eqref{1.1} and \eqref{1.3}. Finally, in Section 4 we prove Theorem \ref{the-1}.

Throughout this paper, we will use $C$ to denote a positive constant, whose value may change from line to line, and write $C_{n,\alpha}$ or $C(n,\alpha)$ to emphasize the dependence of a constant on $n$ and $\alpha$.  For $p\in[1,\infty]$, we denote $L^p(\R^n)$ the standard $L^p$-space and its  norm by  $\|\cdot\|_{L^p(\R^n)}$.
For $s\geq0$, we use notations $\dot{H}^s(\R^n)$ and $H^s(\R^n)$ to denote the homogeneous and nonhomogeneous Sobolev space of $s$ order,
whose endowed norms are denoted by $\|\cdot\|_{\dot{H}^s(\R^n)}=\|\Lambda^s(\cdot)\|_{L^2(\R^n)}$ and $\|\cdot\|_{H^s(\R^n)}=\|\cdot\|_{L^2(\R^n)}+\|\Lambda^s(\cdot)\|_{L^2(\R^n)}$, respectively (see \cite{[BCD]} for more details). Here the fractional Laplacian $\Lambda^s=(-\Delta)^{\frac{s}{2}}$ with $s\in \R$ is defined by \eqref{1.2} or  another equivalent definition as follows: For $0<s<2$ and $\varphi\in C^\infty(\R^n)$
\begin{eqnarray}\label{1.8}
(\Lambda^s\varphi)(x)=C_{n,s}P.V.\int_{\R^n}\frac{\varphi(x)-\varphi(y)}{|x-y|^{n+s}}dy,
\end{eqnarray}
where $C_{n,s}$ is a normalization constant (see \cite{[CC]} and \cite{[DPG]}).

\section{Preliminaries}

 The local well-posedness of \eqref{1.1} in the Sobolev space $H^s(\R^n)$ for some appropriate $s>0$ can be found in \cite{[DC]}, which can be stated as
\begin{lemma}
{\rm(i)} Let $\frac{1}{2}<\alpha<1$ and $s>\frac{n}{2}+2$. Then for each $\theta_0\in H^s(\R^n)$, there exists a $T=T(\|\theta_0\|_{H^s(\R^n)})>0$ such that the equation \eqref{1.1} has a unique solution $\theta$ in $C([0,T); H^s(\R^n))\cap {Lip}((0,T);H^{s-1}(\R^n)).$

{\rm(ii)} Let $0<\alpha\leq\frac{1}{2}$ and $s>\frac{n}{2}+1$. Then for each $\theta_0\in H^s(\R^n)$, there exists a  $T=T(\|\theta_0\|_{H^s(\R^n)})>0$ such that the equation \eqref{1.1} has a unique solution $\theta$ in $C([0,T); H^s(\R^n))\cap {\rm Lip}((0,T);H^{s-1}(\R^n))$. Furthermore, if $T^\ast$ is the first time the solution cannot be continued in $C([0,T^\ast); H^s(\R^n))$, then there necessarily holds $$\int_0^{T^\ast}\|(\mathcal{R}\otimes\mathcal{R})\Lambda^{2\alpha}\theta(\cdot,t)\|_{L^\infty(\R^n)}dt=\infty,$$
\end{lemma}
where $\Big((\mathcal{R}\otimes\mathcal{R})(f)\Big)_{jk}=\mathcal{R}_j\mathcal{R}_k(f)$ denotes the tensor product of the Riesz operators.

The following is  a technical lemma  needed later (see \cite{[V]}).
\begin{lemma}\label{recurrence}
For $C>1$ and $\beta>1$, there exists a constant $C^\ast_0$ such that for every sequence $\{W_k\}^\infty_{k=0}$ satisfying $0<W_0<C^\ast_0$ and for every $k$, $$0\leq W_{k+1}\leq C^kW^\beta_k,$$ we have $\displaystyle\lim_{k\rightarrow\infty}W_k=0.$
\end{lemma}

\begin{remark}\label{2.1}
We can choose the constant $C^\ast_0=C^{-\frac{1}{(\beta-1)^2}}$ in the Lemma \ref{recurrence}.
\end{remark}
The following is an interpolation inequality, which is a direct consequence of the Gagliardo-Nirenberg inequality.
For the convenience of readers, we give a simple proof of this interpolation by the method of the Fourier splitting.
\begin{lemma}\label{interpolation}
Let $\alpha>0$ and $f:\R^n\rightarrow\R$ be a sufficiently regular function. Then there exists a constant $C_{n,\alpha}$ such that
$$\|f\|_{L^2(\R^n)}\leq C_{n,\alpha}\|f\|^{\frac{2\alpha}{n+2\alpha}}_{L^1(\R^n)}\|f\|^{\frac{n}{n+2\alpha}}_{\dot{H}^\alpha(\R^n)}.$$
\end{lemma}
\textbf{Proof}. By virtue of the Plancherel identity of the Fourier transform, for $A>0$ to be determined later, we have
\begin{eqnarray*}\label{2.2}
\|f\|^2_{L^2(\R^n)}
&=&\|\widehat{f}\|^2_{L^2(\R^n)}\\
&=&\int_{|\xi|\leq A}|\widehat{f}(\xi)|^2d\xi+\int_{|\xi|>A}|\widehat{f}(\xi)|^2d\xi\\
&\leq&\omega_n\|f\|^2_{L^1(\R^n)}A^n+A^{-2\alpha}\int_{|\xi|>A}|\xi|^{2\alpha}|\widehat{f}(\xi)|^2d\xi\\
&\leq&\omega_n\|f\|^2_{L^1(\R^n)}A^n+A^{-2\alpha}\|f\|^2_{\dot{H}^\alpha(\R^n)},
\end{eqnarray*}
where $\omega_n=\frac{2\pi^{\frac{n}{2}}}{n\Gamma(\frac{n}{2})}$ is the volume of the unit ball in $\R^n$.
 Letting $A=\Big(\frac {\|f\|_{\dot{H}^\alpha(\R^n)}}{\omega_n^{\frac{1}{2}}\|f\|_{L^1(\R^n)}}\Big)^{\frac{2}{n+2\alpha}}$ yields
$$\|f\|_{L^2(\R^n)}\leq C_{n,\alpha}
\|f\|^{\frac{2\alpha}{n+2\alpha}}_{L^1(\R^n)}
\|f\|^{\frac{n}{n+2\alpha}}_{\dot{H}^\alpha(\R^n)},$$
which concludes the proof of Lemma \ref{interpolation} with $C_{n,\alpha}=2^{\frac{1}{2}}\omega_n^{\frac{\alpha}{n+2\alpha}}.$
\section{Basic properties of solutions}
In this section, we give some basic properties of the smooth solutions to \eqref{1.1} and \eqref{1.3}.
The next lemma is the scaling invariance of the solution to \eqref{1.1}.
\begin{lemma}\label{scaling}
If $\theta$ is a smooth solution to \eqref{1.1}, then for any $\lambda>0$ and $\mu>0$, the rescaled function $$\theta_{\lambda,\mu}(x,t)=(\lambda^{-2\alpha}\mu)\theta(\lambda x,\mu t)$$ is also a solution to \eqref{1.1}.
\end{lemma}
\textbf{Proof}. The proof follows from direct computations.
Firstly,
\begin{eqnarray}\label{3.1}
(\partial_t\theta_{\lambda,\mu})(x,t)=\partial_t((\lambda^{-2\alpha}\mu)\theta(\lambda x,\mu t))=(\lambda^{-2\alpha}\mu^2)(\partial_t\theta)(\lambda x,\mu t).
\end{eqnarray}
Secondly,
\begin{eqnarray}\label{3.2}
(\nabla_x\theta_{\lambda,\mu})(x,t)=\nabla_x((\lambda^{-2\alpha}\mu)\theta(\lambda x,\mu t))=(\lambda^{1-2\alpha}\mu)(\nabla\theta)(\lambda x,\mu t).
\end{eqnarray}
Finally, in view of \eqref{1.4}, we have
\begin{eqnarray}\label{3.3}
u_{\lambda,\mu}(x,t)&=&C_{n,\alpha}P.V.\int_{\R^{n}}\frac{x-y}{|x-y|^{n+2\alpha}}\theta_{\lambda,\mu}(y,t)dy \nonumber \\
&=&(\lambda^{-2\alpha}\mu) C_{n,\alpha}P.V.\int_{\R^{n}}\frac{x-y}{|x-y|^{n+2\alpha}}\theta(\lambda y,\mu t)dy\nonumber\\
&=&(\lambda^{-2\alpha}\mu) C_{n,\alpha}P.V.\int_{\R^{n}}\frac{x-\lambda^{-1}z}{|x-\lambda^{-1}z|^{n+2\alpha}}\theta(z,\mu t)\lambda^{-n}dz \nonumber \\
&=&(\lambda^{-1}\mu) C_{n,\alpha}P.V.\int_{\R^{n}}\frac{(\lambda x-z)}{|\lambda x-z|^{n+2\alpha}}\theta(z,\mu t)dz \nonumber \\
&=&(\lambda^{-1}\mu) u(\lambda x,\mu t),
\end{eqnarray}
where we have made the change of variables $z=\lambda y$ in the above equalities.
Combining \eqref{3.1}, \eqref{3.2} and \eqref{3.3}, we arrive at
\begin{eqnarray*}
&&(\partial_t\theta_{\lambda,\mu}-u_{\lambda,\mu}\cdot\nabla\theta_{\lambda,\mu})(x,t)\\
&=&(\lambda^{-2\alpha}\mu^2)(\partial_t\theta)(\lambda x,\mu t)-(\lambda^{-1}\mu) u(\lambda x,\mu t)\cdot(\lambda^{1-2\alpha}\mu)(\nabla\theta)(\lambda x,\mu t)\\
&=&\lambda^{-2\alpha}\mu^2(\partial_t\theta-u\cdot\nabla\theta)(\lambda x,\mu t)\\
&=&0.
\end{eqnarray*}
The proof of Lemma \ref{scaling} is finished.
\begin{remark}
The solutions of \eqref{1.3} also admit the similar scaling invariance, which involves only one parameter. That is, if $\theta$ is a solution to \eqref{1.3}, then for any $\lambda>0$, the rescaled function $\theta_{\lambda}(x,t)=\lambda^{\gamma-2\alpha}\theta(\lambda x, \lambda^\gamma t)$ is also a solution.
\end{remark}
The following lemma is a maximum principle. We prove it in the spirit of \cite{[CCCF]} and \cite{[CC]}.
\begin{lemma}\label{maximum}
If $\theta$ is a smooth solution to \eqref{1.3} with $\kappa\geq0$, then

{\rm(1)} $M(t)=\displaystyle\sup_{x\in\R^n}\theta(x,t)$ is a Lipchitz continuous and non-increasing function,

{\rm(2)} $m(t)=\displaystyle\inf_{x\in\R^n}\theta(x,t)$ is a Lipchitz continuous and non-decreasing function.
\end{lemma}
\textbf{Proof}. By Sobolev embedding $H^s(\R^n)\hookrightarrow C_0(\R^n)$, $s>\frac{n}{2}$, where $C_0(\R^n)$ is the space of continuous functions on $\R^n$ vanishing at infinity. Therefore if $\theta(\cdot,t)\in C_0(\R^n)$, there exists some point $x_t\in\R^n$ such that $M(t)=\displaystyle\max_{x\in\R^n}\theta(x,t)=\theta(x_t,t).$
For $t_1,t_2\in[0,T)$, by the mean-value theorem, we have
\begin{eqnarray}\label{3.4}
M(t_1)-M(t_2)&=&\theta(x_{t_1},t_1)-\theta(x_{t_2},t_2) \nonumber\\
&=&\theta(x_{t_1},t_1)-\theta(x_{t_1},t_2)+\theta(x_{t_1},t_2)-\theta(x_{t_2},t_2) \nonumber\\
&\leq&\theta(x_{t_1},t_1)-\theta(x_{t_1},t_2) \nonumber\\
&=&(\partial_t\theta)(x_{t_1},\tilde{t})(t_1-t_2)\nonumber\\
&\leq&\|\partial_t\theta\|_{L^\infty(\R^n\times[0,T))}|t_1-t_2|,
\end{eqnarray}
where $\tilde{t}$ is between $t_1$ and $t_2$. Performing the similar argument, we can also obtain
\begin{eqnarray}\label{3.5}
M(t_2)-M(t_1)\leq\|\partial_t\theta\|_{L^\infty(\R^n\times[0,T))}|t_1-t_2|.
\end{eqnarray}
Combining \eqref{3.4} and \eqref{3.5}, we arrive at
\begin{eqnarray*}
|M(t_1)-M(t_2)|\leq\|\partial_t\theta\|_{L^\infty(\R^n\times[0,T))}|t_1-t_2|,
\end{eqnarray*}
which exactly shows that $M(t)$ is Lipschitz continuous in time.
By H. Rademacher's theorem (see \cite{[EG]}), $M(t)$ is differentiable almost everywhere.
Fix $t$ such that $M$ is differentiable at $t$. Under the hypothesis of regularity of $\theta$, we may choose $x_t$ such that $M(t)=\theta(x_t,t).$  By the compactness argument, there exist a positive sequence $\{t_j\}$  and a $x_\ast=x_\ast(t)\in\R^n$ such that $x_{t+t_j}\rightarrow x_\ast$ in $\R^n$ as $t_j\rightarrow0$.
Therefore, it is clear that we have $M(t)=\theta(x_\ast,t)$.
To conclude the proof, we finally consider the finite difference
\begin{eqnarray*}
&&\frac{M(t+t_j)-M(t)}{t_j}\\
&=&\frac{\theta(x_{t+t_j},t+t_j)-\theta(x_\ast,t)}{t_j}\\
&=&\frac{\theta(x_{t+t_j},t+t_j)-\theta(x_{t+t_j},t)}{t_j}+\frac{\theta(x_{t+t_j},t)-\theta(x_\ast,t)}{t_j}\\
&\leq&\frac{\theta(x_{t+t_j},t+t_j)-\theta(x_{t+t_j},t)}{t_j}\\
&=&(\partial_t\theta)(x_{t+t_j},t+\delta t_j)\\
&=&-u(x_{t+t_j},t+\delta t_j)\cdot(\nabla\theta)(x_{t+t_j},t+\delta t_j)-\kappa(\Lambda^\gamma\theta)(x_{t+t_j},t+\delta t_j),
\end{eqnarray*}
for some $\delta\in(0,1)$ and we have used the mean-value theorem.
Taking a limit with $t_j\rightarrow0$ leads to the following inequality
\begin{eqnarray}\label{3.6}
M^\prime(t)&\leq&-u(x_\ast,t)\cdot(\nabla\theta)(x_\ast,t)-\kappa(\Lambda^\gamma\theta)(x_\ast,t) \nonumber \\
&=&-\kappa(\Lambda^\gamma\theta)(x_\ast,t) \nonumber\\
&=&-\kappa C_{n,\gamma}P.V.\int_{\R^n}\frac{\theta(x_\ast,t)-\theta(y,t)}{|x_\ast-y|^{n+\gamma}}dy \nonumber \\
&\leq&0,
\end{eqnarray}
where we have used \eqref{1.8}.
Since \eqref{3.6} holds true at almost every time $t$, we may conclude that $M(t)$ is non-increasing.
By a completely analogous argument, we can obtain that $m(t)$ is also Lipschitz and non-decreasing.
The proof of Lemma \ref{maximum} is complete.

As mentioned above, the velocity field $u$ in \eqref{1.1} is not divergence-free in general. This property makes the integral of the solutions is dissipative, which is
\begin{lemma}\label{identity}
If $\theta$ is a smooth solution to \eqref{1.1}, then $$\frac{d}{dt}\int_{\R^{n}}\theta(x,t)dx=-\|\theta(\cdot,t)\|^2_{\dot{H}^\alpha(\R^n)}.$$
\end{lemma}
\textbf{Proof}. Integrating the first equation in \eqref{1.1} on $\R^n$, integrating by parts and utilizing $\Delta=-\Lambda^2$, we obtain
\begin{eqnarray*}
\frac{d}{dt}\int_{\R^{n}}\theta(x,t)dx&=&-\int_{\R^{n}}u(x,t)\cdot\nabla\theta(x,t)dx\\
&=&\int_{\R^{n}}\theta(x,t) ({\rm div} u)(x,t)dx\\
&=&\int_{\R^{n}}\theta {\rm div} (\nabla\Lambda^{-2+2\alpha}\theta)dx\\
&=&\int_{\R^{n}}\theta(x,t)\triangle\Lambda^{-2+2\alpha}\theta(x,t)dx\\
&=&-\int_{\R^{n}}\theta(x,t)\Lambda^2\Lambda^{-2+2\alpha}\theta(x,t)dx\\
&=&-\int_{\R^{n}}\theta(x,t)\Lambda^{2\alpha}(\theta)(x,t)dx\\
&=&-\|\theta(\cdot,t)\|^2_{\dot{H}^\alpha(\R^n)},
\end{eqnarray*}
which concludes the proof of Lemma \ref{identity}.
\section{Proof of Theorem \ref{the-1}}
In this section, we prove Theorem \ref{the-1}. Before that, we prove a lemma via the De Giorgi iteration technique.
\begin{lemma}\label{control}
If $\theta$ is a smooth solution to \eqref{1.1},  then there exists a constant $\varepsilon_0=\varepsilon_0(n,\alpha)>0$ such that if $\int_{\R^n}\theta(x,0)^+dx\leq\varepsilon_0$ then $\theta(x,1)\leq1$ for all $x\in\R^n.$
\end{lemma}
\textbf{Proof}. Define
\begin{eqnarray}\label{4.1}
\theta_k(x,t):=\Big(\theta(x,t)-C_k\Big)^+,
\end{eqnarray}
where the level value $C_k:=1-2^{-k}$ for $k=0, 1, 2, \cdots $. Here $\Big(\theta(x,t)-C_k\Big)^+$ denotes the positive part of the function $\theta-C_k$.
Clearly, for any $k\geq0$, the truncation $\theta_k$ satisfies the similar transport equation as $\theta$, that is,
\begin{eqnarray}\label{4.2}
\partial_t\theta_k+u\cdot\nabla\theta_k=0.
\end{eqnarray}
Integrating \eqref{4.2} on $\R^n$ and making integration by parts, we have
\begin{eqnarray}\label{4.3}
\frac{d}{dt}\int_{\R^{n}}\theta_k(x,t)dx
&=&-\int_{\R^{n}}u(x,t)\cdot\nabla\theta_k(x,t)dx \nonumber \\
&=&\int_{\R^{n}}\theta_k(x,t)({\rm div}u)(x,t)dx \nonumber \\
&=&-\int_{\R^{n}}\theta_k(x,t)\Lambda^{2\alpha}\theta(x,t)dx \nonumber \\
&=&-\int_{\R^{n}}\theta_k(x,t)\Lambda^{2\alpha}\theta_k(x,t)dx
+\int_{\R^{n}}\theta_k(x,t)\Lambda^{2\alpha}(\theta_k-\theta)(x,t)dx \nonumber \\
&=&-\int_{\R^{n}}\theta_k(x,t)\Lambda^{2\alpha}\theta_k(x,t)dx
+\int_{\R^{n}}\theta_k(x,t)\Lambda^{2\alpha}(\theta_k-\theta+C_k)(x,t)dx \nonumber \\
&=&-\|\theta_k(\cdot,t)\|^2_{\dot{H}^\alpha(\R^n)}+\int_{\R^{n}}\theta_k(x,t)\Lambda^{2\alpha}(\theta_k-\theta+C_k)(x,t)dx,
\end{eqnarray}
where in \eqref{4.3} we have used the fundamental fact $$\Lambda^{2\alpha}(C_k)(x)\equiv0,$$
which is a direct consequence of \eqref{1.8}.
Use again \eqref{1.8} and \eqref{4.1} to obtain
\begin{eqnarray}\label{4.4}
&&\int_{\R^{n}}\theta_k(x,t)\Lambda^{2\alpha}(\theta_k-\theta+C_k)(x,t)dx \nonumber\\
&=&\int_{\{\theta_k(x,t)>0\}}\theta_k(x,t)\Lambda^{2\alpha}(\theta_k-\theta+C_k)(x,t)dx \nonumber\\
&=&\int_{\{\theta(x,t)>C_k\}}\theta_k(x,t)\Big(C_{n,\alpha}\int_{\R^{n}}\frac{[\theta_k-\theta+C_k](x,t)-[\theta_k-\theta+C_k](y,t)}{|x-y|^{n+2\alpha}}dy\Big)dx\nonumber\\
&=&\int_{\{\theta(x,t)>C_k\}}\theta_k(x,t)\Big(C_{n,\alpha}\int_{\R^{n}}\frac{-[\theta_k(y,t)-\theta(y,t)+C_k]}{|x-y|^{n+2\alpha}}dy\Big)dx\nonumber\\
&=&\int_{\{\theta(x,t)>C_k\}}\theta_k(x,t)\Big(C_{n,\alpha}\int_{\{\theta(y,t)\leq C_k\}}\frac{-[\theta_k(y,t)-\theta(y,t)+C_k]}{|x-y|^{n+2\alpha}}dy\Big)dx\nonumber\\
&=&\int_{\{\theta(x,t)>C_k\}}\theta_k(x,t)\Big(C_{n,\alpha}\int_{\{\theta(y,t)\leq C_k\}}\frac{\theta(y,t)-C_k}{|x-y|^{n+2\alpha}}dy\Big)dx\nonumber\\
&\leq&0.
\end{eqnarray}
It follows from \eqref{4.3} and \eqref{4.4} that
\begin{eqnarray}\label{4.5}
\frac{d}{dt}\int_{\R^{n}}\theta_k(x,t)dx\leq-\|\theta_k(\cdot,t)\|^2_{\dot{H}^\alpha(\R^n)}.
\end{eqnarray}
for  $k=0, 1, 2, \cdots$.


To adapt to the De Giorgi iteration technique, we will construct a sequence of time $\{t_k\}^\infty_{k=0}$ by the induction argument such that for all $k>0$, $t_k\in(t_{k-1},C_k)$ and
\begin{eqnarray}\label{4.6}
\|\theta_{k-1}(\cdot,t_k)\|^2_{\dot{H}^\alpha(\R^n)}\leq2^kW_{k-1},
\end{eqnarray}
with $W_k:=\|\theta_k(\cdot,t_k)\|_{L^1(\R^n)}$ and $C_k$ defined as above.

To this end, we proceed in two steps as follows.

\textbf{Step 1.} Let $t_0=0$. Integrating \eqref{4.5} with $k=0$ directly on $[0,\frac{1}{2}]$, we obtain
\begin{eqnarray}\label{4.7}
\int^\frac{1}{2}_0\|\theta_0(\cdot,t)\|^2_{\dot{H}^\alpha(\R^n)}dt&\leq&\int_{\R^n}\theta_0(x,0)dx-\int_{\R^n}\theta_0(x,\frac{1}{2})dx\nonumber\\
&\leq&\int_{\R^n}\theta_0(x,0)dx\nonumber\\
&=&W_0.
\end{eqnarray}
On the other hand, by the mean value theorem, there exists $t_1\in(0,\frac{1}{2})$ such that
\begin{eqnarray}\label{4.8}
\int^\frac{1}{2}_0\|\theta_0(\cdot,t)\|^2_{\dot{H}^\alpha(\R^n)}dt=\frac{1}{2}\|\theta_0(\cdot,t_1)\|^2_{\dot{H}^\alpha(\R^n)}.
\end{eqnarray}

Consequently, in view of \eqref{4.7} and \eqref{4.8}, we obtain $\|\theta_0(\cdot,t_1)\|^2_{\dot{H}^\alpha(\R^n)}\leq2W_0.$ This is exactly the construction of $t_1$.

\textbf{Step 2.}
Assume that we have constructed $t_k\in(t_{k-1},C_k)$ up to some value $k\geq1$.
Then, integrating \eqref{4.5} on the interval $[t_k,C_{k+1}]$ , we can obtain
\begin{eqnarray}\label{4.9}
\int^{C_{k+1}}_{t_k}\|\theta_k(\cdot,t)\|^2_{\dot{H}^\alpha(\R^n)}dt&\leq&\int_{\R^n}\theta_k(x,t_k)dx-\int_{\R^n}\theta_k(x,C_{k+1})dx\nonumber\\
&\leq&\int_{\R^n}\theta_k(x,t_k)dx\nonumber\\
&=&W_k.
\end{eqnarray}
On the other hand, by the mean value theorem, there exists $t_{k+1}\in(t_k,C_{k+1})$ such that
\begin{eqnarray}\label{4.10}
\int^{C_{k+1}}_{t_k}\|\theta_k(\cdot,t)\|^2_{\dot{H}^\alpha(\R^n)}dt
=(C_{k+1}-t_k)\|\theta_k(\cdot,t_{k+1})\|^2_{\dot{H}^\alpha(\R^n)}.
\end{eqnarray}
Combining \eqref{4.9} and \eqref{4.10}, we arrive at
\begin{eqnarray}\label{4.11}
\|\theta_k(\cdot,t_{k+1})\|^2_{\dot{H}^\alpha(\R^n)}&\leq&\frac{W_k}{C_{k+1}-t_k} \nonumber\\
&\leq&\frac{W_k}{C_{k+1}-C_k}\\
&=&2^{k+1}W_k, \nonumber
\end{eqnarray}
where in \eqref{4.11} we have used the induction assumption $t_k\in(t_{k-1},C_k).$
Therefore, we have obtained the desired sequence $\{t_k\}^\infty_{k=0}$.

Recalling \eqref{4.1} and the definition of $W_k$, we have
\begin{eqnarray}\label{4.12}
W_{k+1}&=&\int_{\R^n}\theta_{k+1}(x,t_{k+1})dx \nonumber \\
&=&\int_{\R^n}\Big(\theta(x,t_{k+1})-C_{k+1}\Big)^+dx \nonumber \\
&=&\int_{\R^n}\Big(\theta(x,t_{k+1})-C_k-2^{-k-1}\Big)^+dx \nonumber \\
&=&\int_{\R^n}\Big(\Big(\theta(x,t_{k+1})-C_k\Big)^+-2^{-k-1}\Big)^+dx \nonumber \\
&=&\int_{\R^n}\Big(\theta¡ª_k(x,t_{k+1})-2^{-k-1}\Big)^+dx \nonumber \\
&=&\int_{\{\theta_k(x,t_{k+1})>2^{-k-1}\}}\Big(\theta¡ª_k(x,t_{k+1})-2^{-k-1}\Big)dx \nonumber \\
&\leq&\int_{\{\theta_k(x,t_{k+1})>2^{-k-1}\}}\theta¡ª_k(x,t_{k+1})dx\nonumber\\
&\leq&\|\theta_k(\cdot,t_{k+1})\|_{L^2(\R^n)}|\{x\in\R^{n}:\theta¡ª_k(x,t_{k+1})>2^{-k-1}\}|^\frac{1}{2}\nonumber\\
&\leq&2^{k+1}\|\theta_k(\cdot,t_{k+1})\|^2_{L^2(\R^n)},
\end{eqnarray}
where we have used Cauchy-Schwartz and Chebechev inequality in the last two inequalities.
Integrating \eqref{4.5} again on the time interval $[t_k,t_{k+1}]$, we can obtain
\begin{eqnarray}\label{4.13}
\|\theta¡ª_k(\cdot,t_{k+1})\|_{L^1(\R^n)}\leq\|\theta¡ª_k(\cdot,t_k)\|_{L^1(\R^n)}.
\end{eqnarray}
It follows from Lemma \ref{interpolation}, \eqref{4.12}, \eqref{4.13} and \eqref{4.6} that
\begin{eqnarray}\label{4.14}
W_{k+1}&\leq&2^{k+1}\|\theta_k(\cdot,t_{k+1})\|^2_{L^2(\R^n)}\nonumber\\
&\leq&2^{k+1}C_{n,\alpha}\|\theta_k(\cdot,t_{k+1})\|^{\frac{4\alpha}{n+2\alpha}}_{L^1(\R^n)}
\|\theta_k(\cdot,t_{k+1})\|^{\frac{2n}{n+2\alpha}}_{\dot{H}^\alpha(\R^n)}\nonumber\\
&\leq&2^{k+1}C_{n,\alpha}\|\theta_k(\cdot,t_k)\|^{\frac{4\alpha}{n+2\alpha}}_{L^1(\R^n)}
\|\theta_k(\cdot,t_{k+1})\|^{\frac{2n}{n+2\alpha}}_{\dot{H}^\alpha(\R^n)}\nonumber\\
&\leq&C_{n,\alpha}2^{k+1}W_k^{\frac{4\alpha}{n+2\alpha}}(2^{k+1}W_k)^{\frac{n}{n+2\alpha}}\nonumber\\
&\leq&C_{n,\alpha}2^{\frac{2n+2\alpha}{n+2\alpha}k}W_k^{\frac{n+4\alpha}{n+2\alpha}}.
\end{eqnarray}
With the aid of Lemma \ref{recurrence}, Remark \ref{2.1} and \eqref{4.14}, we can choose $\varepsilon_0=\frac{2^{-\frac{(n+\alpha)(n+2\alpha)}{2\alpha^2}}}{C_{n,\alpha}^{\frac{n+2\alpha}{2\alpha}}}$.
Furthermoer, if $\int_{\R^n}\theta(x,0)^+dx\leq\varepsilon_0$, that is $W_0\leq\varepsilon_0$, then $\displaystyle\lim_{k\rightarrow\infty}W_k=0$.
This implies that $\theta(x,t_{\infty})\leq1$, where $t_\infty:=\displaystyle\lim_{k\rightarrow\infty}t_k\leq1$.
By Lemma \ref{maximum}, we obtain $\theta(x,1)\leq1$ for all $x\in\R^n.$
The proof of Lemma \ref{control} is finished.


Based on Lemma \ref{control}, we have
\begin{lemma}\label{decay}
Suppose that $\theta$ is a smooth solution to \eqref{1.1} with a initial data $\theta_0$ satisfying its positive part $\theta_0^+$ in $L^1(\R^n)$. Then there exists a constant $\varepsilon_0>0$ only depending on $n$ and $\alpha$, such that for every $T>0$
we have
\begin{eqnarray*}
\sup_{x\in\R^n}\theta(x,T)\leq \Big(\frac{\|\theta_0^+\|_{L^1(\R^n)}}{\varepsilon_0T^{\frac{n}{2\alpha}}}\Big)^{\frac{2\alpha}{n+2\alpha}}.
\end{eqnarray*}
\end{lemma}
\textbf{Proof}. Let $\varepsilon_0$ be the constant in Lemma \ref{control}.
Consider the rescaled function $\omega$, defined by $$\omega(x,t)=(\lambda^{-2\alpha}\mu)\theta(\lambda x,\mu t),$$
with $\lambda=\Big(\frac{T\|\theta_0^+\|_{L^1(\R^n)}}{\varepsilon_0}\Big)^{\frac{1}{n+2\alpha}}$ and $\mu=T.$ Then $\omega(x,1)=(\lambda^{-2\alpha}\mu)\theta(\lambda x,T)$,
we further have
\begin{eqnarray}\label{4.15}
\theta(x,T)=\frac{\lambda^{2\alpha}}{\mu}\omega(\frac{x}{\lambda},1).
\end{eqnarray}
From Lemma \ref{scaling}, $\omega$ is also a solution to \eqref{1.1} with the initial data $\omega(x,0)=(\lambda^{-2\alpha}\mu)\theta_0(\lambda x)$.
Note that $\int_{\R^n}\omega(x,0)^+dx=\lambda^{-n-2\alpha}\mu\|\theta_0^+\|_{L^1(\R^n)}=\varepsilon_0$.
Then, applying Lemma \ref{control}, we obtain $\displaystyle\sup_{x\in\R^n}\omega(x,1)\leq1$.
Substituting this bound into \eqref{4.15}, we obtain
\begin{eqnarray*}
\sup_{x\in\R^n}\theta(x,T)
&\leq&\frac{1}{T}\Big(\frac{T\|\theta_0^+\|_{L^1(\R^n)}}{\varepsilon_0}\Big)^{\frac{2\alpha}{n+2\alpha}}\\
&=&\Big(\frac{\|\theta_0^+\|_{L^1(\R^n)}}{\varepsilon_0T^{\frac{n}{2\alpha}}}\Big)^{\frac{2\alpha}{n+2\alpha}},
\end{eqnarray*}
which is the desired decay estimate.

With help of the above lemmas, we are now ready to prove Theorem \ref{the-1}.

\textbf{Proof of Theorem \ref{the-1}.} We will argue by contradiction. Suppose that $\theta$ is a smooth global solution to \eqref{1.1}. Then the function $\displaystyle\sup_{x\in\R^n}\theta(x,t)$ must be constant.
Indeed, pick a point $x_M\in \R^n$ such that $\theta_0$ attains its maximum at $x_M$ and let it evolve according to the following ordinary differential equation
\begin{equation*}
\left\{\ba
&\frac{d}{dt}x(t)=u(x(t),t), \\
&x(t)|_{t=0}=x_M.\ea\ \right.
\end{equation*}
Then, it deduces that
\begin{eqnarray*}
\frac{d}{dt}\theta(x(t),t)&=&(\partial_t\theta)(x(t),t)+\nabla\theta(x(t),t)\cdot\frac{d}{dt}x(t)\\
&=&(\partial_t\theta)(x(t),t)+u(x(t),t)\cdot\nabla\theta(x(t),t)\\
&=&(\partial_t\theta+u\cdot\nabla\theta)(x(t),t)\\
&=&0,
\end{eqnarray*}
which shows that $\theta(x,t)$ keeps constant on the trajectory line, that is,
$$\theta(x(t),t)=\theta(x(0),0)=\theta_0(x_M)=\sup_{x\in\R^n}\theta_0(x).$$
Furthermore,
\begin{eqnarray}\label{4.16}
\sup_{x\in\R^n}\theta(x,t)\geq\theta(x(t),t)=\sup_{x\in\R^n}\theta_0(x)\geq\sup_{x\in\R^n}\theta(x,t),
\end{eqnarray}
where the last inequality in \eqref{4.16} follows from Lemma \ref{maximum}.
It deduces that $\displaystyle\sup_{x\in\R^n}\theta(x,t)$ is the constant function $\displaystyle\sup_{x\in\R^n}\theta_0(x)$ independent of $t$.
However, according to Lemma \ref{decay}, we obtain, for all time $t>0$,
\begin{eqnarray}\label{4.17}
0<\sup_{x\in\R^n}\theta_0(x)=\sup_{x\in\R^n}\theta(x,t)\leq \Big(\frac{\|\theta_0^+\|_{L^1(\R^n)}}{\varepsilon_0t^{\frac{n}{2\alpha}}}\Big)^{\frac{2\alpha}{n+2\alpha}}.
\end{eqnarray}
Let $t$ be sufficiently large such that
\begin{eqnarray}\label{4.18}
\Big(\frac{\|\theta_0^+\|_{L^1(\R^n)}}{\varepsilon_0t^{\frac{n}{2\alpha}}}\Big)^{\frac{2\alpha}{n+2\alpha}}
<\frac{1}{2}\sup_{x\in\R^n}\theta_0(x).
\end{eqnarray}
Combining \eqref{4.17} and \eqref{4.18}, we obtain an obvious contradiction, which then shows that local smooth solution to \eqref{1.1} must develop singularity in finite time.
The proof of Theorem \ref{the-1} is complete.

{\bf Acknowledgements.}
Q. Jiu was partially supported by the National Natural Science Foundation of
China (NNSFC) (No. 11931010, No. 12061003)
and key research project of the Academy for Multidisciplinary Studies of CNU.


\end{document}